\documentclass[12pt]{amsart}

\usepackage{amsmath}

\def\nbd#1#2{{\mathcal N} (#1;#2)}

\newcommand{\BB}{{\mathbb B}}
\newcommand{\RR}{{\mathbb R}}
\newcommand{\CC}{{\mathbb C}}

\def\Re{\hbox{Re}\,}
\def\Im{\hbox{Im}\,}

\def\O{\Omega}

\def\auto{\hbox{\rm Aut}\,(\Omega)}
\def\bo{\partial \Omega}

\def\HollowBox #1#2{{\dimen0=#1 \advance\dimen0 by -#2
       \dimen1=#1 \advance\dimen1 by #2
        \vrule height #1 depth #2 width #2
        \vrule height 0pt depth #2 width #1
        \llap{\vrule height #1 depth -\dimen0 width \dimen1}%
       \hskip -#2
       \vrule height #1 depth #2 width #2}}
\def\BoxOpTwo{\mathord{\HollowBox{6pt}{.4pt}}\;}

\def\cE{{\mathcal E}}

\def\cH{{\mathcal H}}

\def\origin{\mathbf 0}
\def\ee{{\mathbf e}}
\def\ff{{\mathbf f}}
\def\pp{{\mathbf p}}
\def\qq{{\mathbf q}}
\def\vv{{\mathbf v}}
\def\ww{{\mathbf w}}
\def\xx{{\mathbf x}}
\def\yy{{\mathbf y}}
\def\zz{{\mathbf z}}

\def\Dbar{\overline{D}}

\newtheorem{theorem}{Theorem}[section]
\newtheorem{lemma}[theorem]{Lemma}
\newtheorem{proposition}[theorem]{Proposition}
\newtheorem{corollary}[theorem]{Corollary}

\newtheorem{definition}[theorem]{DEFINITION}

\title[Characterization of the Hilbert ball]
{Characterization of the Hilbert ball by its automorphism group}

\author{Kang-Tae Kim \and Steven G. Krantz}

\begin{document}
\maketitle

\begin{quote}
\sl Let $\Omega$ be a bounded, convex domain in a separable
Hilbert space.  The authors prove a version of the theorem of Bun
Wong, which asserts that if such a domain admits an
automorphism orbit accumulating at a strongly pseudoconvex
boundary point then it is biholomorphic to
the ball.  Key ingredients in the proof are a new localization
argument using holomorphic peaking functions and
the use of new ``normal families'' arguments in the construction
of the limit biholomorphism.
\end{quote}

\section{Introduction}

Bun Wong [WON] proved in 1977 that a strongly pseudoconvex domain
in $\CC^n$ with non-compact automorphism group must be the ball.
In the intervening 22 years, there has been considerable work
on this circle of ideas.  Rosay [ROS] showed that (in order
to draw the same conclusion) the domain can
be allowed to be bounded, but arbitrary, and one need only assume that
some boundary orbit accumulation point of the automorphism group action
be strongly pseudoconvex.

Greene and Krantz [GRK1] explored what happens when the orbit
accumulation point is weakly pseudoconvex.  Although the ideas
developed in that paper have now been superseded by more modern methods that
center around the technique of scaling (see [IK1], [Ki4], [KIKR]
for a discursive treatment of the scaling method), certainly that
paper ushered in a new avenue of inquiry.  Subsequent work
by Bedford/Pinchuk [BP1--BP4], by Frankel [Fr],
by Kim [Ki2--Ki6], by Fu/Isaev/Krantz ([IK2--IK4], [FIK1--FIK2]),
by Berteloot and Coeur\'{e}
([Ber1--Ber3], [BeCo]), and by many others has fleshed
out this circle of ideas.

It should be stressed, however, that all work to date has been in
complex spaces of finite dimensions.  The most complete results are
available in $\CC^2$.  The work in dimensions 3 and higher continues
to develop.  Nevertheless, some questions remain open in {\it all}
dimensions 2 and higher.

In the present work we begin to explore a domain in an infinite
dimensional complex Hilbert space; we assume that the domain
possesses an automorphism group orbit accumulating at a boundary
point.  [Notice that the mere non-compactness of the automorphism
group with respect to the compact-open topology does not give the
(rigid) structure that we require.]  In particular, we shall
formulate and prove a version of the Wong/Rosay theorem
in that setting.

We are pleased to thank Yun-Sung Choi, Jaesung Lee, John M${}^{\rm
c}$Carthy, and Richard Rochberg for helpful discussions.
The work of Laszlo Lempert sparked our interest in the holomorphic
function theory of infinitely many variables, and his helpful advice
is also much appreciated.

The first author's research was funded partly by Grant 981-0104-018-2
and Interdisciplinary Research Program Grant 1999-2-102-003-5 of the
Korean Science and Engineering Foundation. The second author was
supported in part by NSF Grants DMS-9531967 and DMS-9631359.

\section{Preliminaries}

Let $\cH$ be an infinite-dimensional, separable Hilbert space over the
complex numbers.

\subsection{Derivatives and $C^k$ smoothness}

First we consider the differentiability of functions.  Let $W
\subseteq \cH$ be an open subset.  Consider $u:W
\to \cE$, where $\cE$ is some complex Hilbert space.
For $\qq \in W$ and $\vv_1, \ldots, \vv_k \in \cH$,
we consider the following derivatives:
\begin{eqnarray*}
du(\qq;\vv_1)
  & = & \lim_{\RR \ni \epsilon \to 0}
\frac{1}{\epsilon} \left( u(\qq+\epsilon \vv_1) - u(\qq) \right) \\
d^2 (\qq;\vv_1, \vv_2)
  & = &
  \lim_{\RR \ni \epsilon \to 0} \frac{1}{\epsilon}
\left( du(\qq+\epsilon \vv_2;\vv_1) - du(\qq;\vv_1) \right)
  \\
  & \vdots & \\
d^k ({\bf q}; {\bf v}_1, \ldots, {\bf v}_k)
   & = & \\
&&  \hspace{-80pt}
  \lim_{\RR \ni \epsilon \to 0}
  \frac{1}{\epsilon} \left( d^{k-1}
  u(\qq+\epsilon \vv_k; \vv_1, \ldots, \vv_{k-1})
    - d^{k-1}u(\qq; \vv_1, \ldots, \vv_{k-1}) \right).
\end{eqnarray*}
We say that the mapping $u$ is {\it $C^k$ smooth} (for $k \ge 1$)
if each $d^\ell u$, $\ell=1,\ldots,k$, exists and is a
continuous map of
$W \times \underbrace{\cH \times \ldots \cH}_{\ell\ \rm copies}$
into $\cE$.

In this case, it is known that each $d^\ell u$
is a symmetric $\ell$-linear map over $\RR$.
\medskip \\

We may also consider the following derivatives:
\begin{eqnarray*}
D u(\qq;\vv) & = & \frac{du(\qq;\vv) - i\, du(\qq; i\vv)}{2} \\
\Dbar u(\qq;\vv) & = & \frac{du(\qq;\vv) + i\, du(\qq; i\vv)}{2}
\end{eqnarray*}
Since a $C^1$ smooth function that is holomorphic in every
complex direction is holomorphic (this is the {\it Gateaux
 definition
of holomorphicity}, for which see [Mu]), we see that every $C^1$ smooth
mapping $u$ with $\Dbar u$ identically zero is holomorphic.
An alternative, and equivalent, definition is to demand that
$u$ be holomorphic on each finite-dimensional slice.  Both of
these are equivalent to specifying that $u$ have a power series
expansion about each point.
\rm
\medskip  \\

\subsection{Strong pseudoconvexity}

By a \it domain \rm in $\cH$, we mean a connected open subset of $\cH$.
If $\O$ is a domain in $\cH$, then we denote its boundary by $\bo$.  We say
that a boundary point $\pp \in \bo$ is $C^k$ {\it smooth} ($k\ge 1$) if $\pp$
admits an open neighborhood $U$ in $\cH$ and a $C^k$ smooth function
$\rho:U \to \RR$ satisfying the following conditions:
\medskip  \\

\begin{itemize}
\item[(\romannumeral1)]
$\O \cap U = \{ \zz \in U \mid \rho (\zz) < 0 \}$,
\item[(\romannumeral2)]
$\bo \cap U = \{ \zz \in U \mid \rho (\zz) = 0 \}$,
\item[(\romannumeral3)]
$U \setminus \overline{\O} = \{ \zz \in U \mid \rho (\zz) > 0 \}$
and
\item[(\romannumeral4)]
the differential $d\rho_\qq$ is a non-zero functional for every
$\qq \in \bo \cap U$.
\end{itemize}
\vspace*{.1in}

\noindent
The function $\rho$ is called a {\it local defining function} for $\O$
at $p$.  The local defining function is not unique in general, but
it is in effect unique since two given defining functions will differ
only multiplicatively by a positive, smooth function.
\medskip \\

Now we would like to define strong pseudoconvexity in the sense
of Levi:

\begin{definition} \rm
Consider the case in which the boundary point $\pp \in \partial \Omega$
under consideration
admits a $C^2$ local defining function $\rho$ defined on an open
neighborhood $U$ of $\pp$.  Then we call the Hermitian
form
$$
\Lambda\rho (\pp;\vv,\ww) \equiv D\Dbar\rho (\pp;\vv,\ww)
$$
the {\it Levi form of} $\rho$ at $\pp$.  We say that a point
$\pp \in \bo$ is a {\it strongly pseudoconvex boundary
point} if there exist a $C^2$ smooth local defining function
$\rho$ for $\Omega$ in a neighborhood of $\pp$ and a constant
$C > 0$ satisfying the condition
$$
\Lambda\rho (\pp;\vv,\vv) \ge C \| \vv \|^2, ~~\forall \vv \in \cH.
$$
\end{definition}

Note that this definition of strong pseudoconvexity is equivalent
to requiring that $\pp \in \bo$ admit an open neighborhood inside
of which $\bo$ can be biholomorphically transformed to a strongly
convex set.  We will specify and use the precise meaning
of this remark in what follows.

\section{Main Theorem}

We now state the principal result of the present paper.

\begin{theorem} \sl
Let $\Omega$ be a bounded convex domain in a separable Hilbert
space $\cH$.  Assume that $\Omega$ admits a boundary point
$\pp \in \bo$ at which
\begin{itemize}
\item[(1)] $\bo$ is $C^2$ smooth and strongly pseudoconvex
in a neighborhood of $\pp$, and
\item[(2)] there exist $\qq \in \O$ and $f_j \in \auto$ ($j=1,2,\ldots$)
such that $f_j (\qq)$ converges to $\pp$ in norm as $j \to \infty$.
\end{itemize}
Then $\O$ is biholomorphic to the unit ball $\BB = \{z \in \cH
\mid \|\zz\| < 1 \}$.
\end{theorem}

Theorem 3.1 is strongly analogous to the classical result of
Rosay [ROS].  Its proof, however, is by no means a direct
transliteration of Rosay's proof.  To point out critical differences,
we would like to mention first that the issues of localization in
infinite dimensions are quite delicate and require new arguments.
In addition, the holomorphic invariants based upon the
classical volume elements are
no longer valid.  The geometry of Hilbert space presents
some unfamiliar subtleties, even in the application of normal
families of holomorphic mappings and of scaling methods.

\section{Localization}

It is a classical result of Cartan (see [NAR] and references
therein) that if a bounded domain $\Omega \subseteq \CC^n$ has
non-compact automorphism group then there are a point
${\bf q} \in \Omega$ and a sequence
of automorphisms $f_j \in \hbox{Aut}(\Omega)$ such that
$f_j({\bf q}) \rightarrow {\bf p}$ for some ${\bf p} \in \partial \Omega$.  The
converse is true as well.  Thus in practice, we think of non-compactness
of the automorphism group, and non-compactness of an orbit, as
interchangeable.

However, in the infinite dimensional case, these two notions are not
equivalent.  For instance, the group of rotations (i.e. Hilbert space
isometries) acts on the unit ball $\BB$ in Hilbert space; it is
non-compact, although it does not generate an orbit accumulating at a
boundary point of $\BB$.  Several considerations (in
light of the function theory and the geometric theory of bounded domains)
lead us to believe that the existence of an automorphism orbit
accumulating at a boundary point is a more essential condition than
mere non-compactness of the automorphism group.

Thus we pose the following condition (which clearly plays a pivotal
role in our formulation of Theorem 3.1):
\bigskip  \\

\begin{quote}
\begin{itemize}
\item[(\dag)] \it
Let $\O$ be a bounded domain in $\cH$ with a boundary point
$\pp$ near which $\bo$ is $C^2$ smooth and strongly pseudoconvex.
We assume that there is
a point $\qq \in \O$ and a sequence of holomorphic
automorphisms $f_j \in \auto$ such that $f_j(\qq) \to \pp$
as $j$ tends to infinity.
\end{itemize}
\end{quote}
\vspace*{.2in}

Crucial to the analysis that we have just outlined is, at
least in the finite-dimensional case of the Wong-Rosay theorem, that one
be able to {\it localize} the calculation of a certain
invariant near {\bf p}.  This invariant, denoted in the literature
by ${\mathcal M}$, is the ratio between the Carath\'{e}odory
and Kobayashi volume elements.  More precisely, one wants to know that
if $U$ is a small neighborhood of {\bf p} in $\CC^n$ then
the calculation of ${\mathcal M}(f_j({\bf q}))$ relative
to $\Omega \cap U$ is essentially the same as the calculation of
${\mathcal M}(f_j({\bf q}))$ relative to $\Omega$.  In the
infinite dimensional case, the ${\mathcal M}$ invariant has no suitable
definition and hence is no longer relevant.  However, the localization
argument which relates arbitrarily large subdomains to a local neighborhood
of the orbit accumulation point is still quite essential.

The localization in the classical setting of finite dimensions involves
delicate normal families arguments (see [KRA1]).
In infinite dimensions it is considerably more subtle.  This difficulty
is exacerbated in part by the non-existence of a compact
exhaustion of the domain.

We next present an effective version of localization in the
infinite dimensional setting.
\smallskip \\

We first need a suitable version of Montel's theorem in infinite
dimensions.  We present one now.

\begin{lemma} \sl
Let $U$ be an open set in $\cH$. Let ${\mathcal F} =
\{f_\alpha\}_{\alpha \in {\mathcal A}}$ be a family of holomorphic
(scalar-valued) functions on $U$.  Assume that there is a constant
$M > 0$ such that $|f_\alpha(z)| \leq M$ for all $\alpha \in
{\mathcal A}$ and all $z \in U$.  Then each subsequence
$\{f_{\alpha_j}\}$ has itself a subsequence $\{f_{\alpha_{j_k}}\}$
that converges uniformly on compact subsets of $U$.
\end{lemma}

\noindent {\bf Proof:}   Fix a compact set $K \subset U$. There is
a positive number $\delta_0$ such that $\|k - u\| \geq 2\delta_0$
for all $k \in K$ and $u \not \in U$.  Now for $k\in K$ and any unit
vector $\eta \in \cH$ we consider the disc ${\mathcal D}
 = \{k + \zeta \eta : \zeta \in \CC, |\zeta| < \delta_0\}$.
Applying the standard one-variable Cauchy estimates
on $\mathcal D$, we find that the directional derivative of any
$f_\alpha$ at $k$ in the direction $\eta$ has a uniform
bound---depending on $\delta_0$ and $M$, but not on $\alpha$.
Therefore the elements of ${\mathcal F}$ satisfy a uniform
Lipschitz estimate
$$
\|f_\alpha(z) - f_\alpha(z + h)\| \leq K \cdot \|h\| ,
$$
where $K$ a constant independent of $\alpha$.  As
a result, the family ${\mathcal F}$ is equicontinuous.  Of course
it is equibounded by hypothesis.

As a result of these considerations, the Ascoli-Arzela theorem
(which is valid on a compact metric space)
applies and the required conclusion follows.
\qed
\smallskip \\

\noindent {\bf Remark:}  The Montel theorem that we have just
formulated and proved will certainly suffice for our purposes. But
it is worth noting that the following (formally) more general
statement is now easily derived:

\begin{quote} \sl
Let $U$ be an open set in $\cH$. Let ${\mathcal F} =
\{f_\alpha\}_{\alpha \in {\mathcal A}}$ be a family of holomorphic
(scalar-valued) functions on $U$.  Assume that for each compact
set $K \subset U$ there is a constant $M_K > 0$ such that
$|f_\alpha(z)| \leq M_K$ for all $\alpha \in {\mathcal A}$ and all
$z \in M_K$.  Then each subsequence $\{f_{\alpha_j}\}$ has itself
a subsequence $\{f_{\alpha_{j_k}}\}$ that converges uniformly on
compact subsets of $U$.
\end{quote}

For a proof, we argue as follows (we are indebted to Laszlo
Lempert for this elegant idea).  We claim that each $P \in U$ has
a neighborhood
$$
\nbd{P}{r} \equiv \{z \in \cH : \|z - P\| < r \}
$$
on which the family ${\mathcal F}$ is uniformly bounded.
If not, then inside each $\nbd{P}{1/n}$ there is a point $x_n$ and
a function $f_{\alpha_n} \in {\mathcal F}$ such that
$|f_{\alpha_n}(x_n)| > n$.  Then
$$
L = \{x_1, x_2, \dots\} \cup \{P\}
$$ is a compact set and ${\mathcal F}$ is unbounded on $K$.
That is a contradiction.  As a result, if $K$ is any compact set
in $U$ then we can find a full neighborhood $W$ of $K$ on which
${\mathcal F}$ is uniformly bounded.  Then Lemma 4.1, as stated,
applies on $W$. \hfill
$\BoxOpTwo$
\vspace*{.13in}

Now we have:

\begin{lemma} \sl
Let $U$ be an open neighborhood of $\pp \in \partial \Omega$
and assume that $\xx \in
\O$ satisfies $\displaystyle{\lim_{j\to\infty} f_j (\xx) =
{\bf p}}$.  Assume that $\Omega$ is strongly pseudoconvex
in a neighborhood of $\pp$.
Then there exist a real number $r_\xx > 0$ and an integer $\nu_\xx
> 0$ such that $f_\nu (\nbd{\xx}{r_\xx}) \subset U \cap \Omega$
for every $\nu > \nu_\xx$.
\end{lemma}

\bf Proof: \rm Seeking a contradiction, we assume the contrary.
This hypothesis in particular implies:
\medskip  \\
\begin{quote}
\begin{itemize}
\item[($\ddag$)] \it There exists a subsequence $f_{j_\nu}$ of the
sequence $f_j$ together with a point sequence
$\xx_\nu$ that converges to $\xx$ in norm as $\nu \to \infty$ such
that $f_{j_\nu} (\xx_\nu) \notin U \cap \O$.
\end{itemize}
\end{quote}
\vspace*{.15in}

We shall denote the subsequence by $f_\nu$ for convenience.
\medskip  \\

Let $\{\ee_1, \ee_2, \ldots\}$ be an orthonormal basis for $\cH$.
In what follows, we think of an element ${\mathbf h} \in \cH$ as having
the form ${\bf h} = \sum_j z_j {\bf e}_j$, and we identify $\bf h$ with its
coordinates $(z_1, z_2, \dots)$.  Now an analog of the finite-dimensional
case allows us to consider a quadratic
holomorphic polynomial mapping $G : \cH \to \cH$ which maps $U$
(an open neighborhood of $\pp$ in $\cH$) biholomorphically onto
the open neighborhood $G(U)$ of $G(\pp)$ so that
the following properties hold:
\begin{itemize}
\item[(1)] $\O_U \equiv
G(U \cap \O)$ is defined by the inequality
$$
\Re z_1 > \psi
\left((\Im z_1) \ee_1 + \sum_{j=2}^\infty z_j \ee_j \right) \, ;
$$
\item[(2)] $G(\pp)=\origin$;
\item[(3)] $\O_U$  is supported by the closed subspace (a hyperplane)
of codimension one defined by the equation $\Re z_1 = 0$;
\item[(4)] the function $\psi:G(U) \to \RR$ is non-negative and
strongly convex in the sense that it is $C^2$ smooth and there
exists a constant $C>0$ satisfying
$$
d^2 \psi (\yy; \vv, \vv) \ge C \| \vv \|^2 \hbox{ for every } \yy \in G(U)
\hbox{ and } \vv \in \cH \, .
$$
\end{itemize}

Notice that we may choose the neighborhood $U$ and the positive
number $\epsilon$ so that the closure of the set
$$
V = \biggl \{
\zz = \sum_{j=1}^\infty z_j \ee_j \in \O_U \mid \Re z_1 < \epsilon
\biggr \}
$$
is contained in $G(U)$.  Now let
$$
\Sigma = \biggl
\{ \zz = \sum_{j=1}^\infty z_j \ee_j \in \overline{\O}_U \mid \Re
z_1 = \epsilon \biggr \}.
$$
Note that $\Sigma$ separates $\O_U$
into two disjoint subsets. Therefore $G^{-1}(\Sigma)$ also
separates $U \cap \O$ into two disjoint subsets.

Now consider the straight line segment $\gamma_\nu$ joining $\xx$
and $\xx_\nu$.  The segment $\gamma_\nu$ is contained in $\O$ for
sufficiently large values of $\nu$, since $\xx_\nu$ converges to
$\xx$ in norm as $\nu \to \infty$.  Then the curve $f_\nu \circ
\gamma_\nu$ connects $f_\nu (\xx)$ and $f_\nu (\xx_\nu)$.  For
sufficiently large values of $\nu$, we have that $f_\nu (\xx) \in
(G|_U)^{-1} (V)$ by hypothesis, and that $f_\nu (\xx_\nu)
\notin (G|_U)^{-1} (V)$ by ($\ddag$).  In particular, the set
$$
\hbox{Image}\,(f_\nu \circ \gamma_\nu) \setminus (G|_U)^{-1}
(\Sigma)
$$
is disconnected.  Therefore we can easily conclude by connectivity
that there exists a point $\tilde \xx_\nu \in \gamma_\nu$ such
that $f_\nu (\tilde \xx_\nu) \in (G|_U)^{-1} (\Sigma)$.  Notice
that $\tilde \xx_\nu$ also converges to $\xx$ in norm.

Consider the sequence $G (f_\nu (\tilde \xx_\nu))$.  This sequence
consists of points in $\Sigma$.  By the Banach-Alaoglu theorem, it
has a subsequence that converges weakly.  By an abuse of notation
let us say that the sequence $G (f_\nu (\tilde \xx_\nu))$
converges to the point $\yy^* \in \cH$ weakly.  Since $\Sigma$ is
closed, bounded, and convex, we have $\yy^* \in \Sigma$.

Now observe that, for an arbitrary $\vv \in \cH$, we have
\begin{eqnarray*}
- \langle \yy^* , \vv \rangle & = & \langle 0 - \yy^*, \vv \rangle
\\ & = & \lim_{\nu\to\infty}
  \langle G \circ f_\nu (\xx) - G \circ f_\nu (\tilde{\xx}_\nu), \vv \rangle
\end{eqnarray*}
On the other hand, notice that the set $ Y = \{\xx\} \cup
\{\xx_\nu \mid \nu = 1,2,\ldots\}$ is compact in the strong
topology, since $\xx_\nu \to \xx$ in norm. Let $\ell_\vv(\zz) \equiv
\langle \zz, {\bf v}\rangle$. Then Lemma 4.1 implies that every
subsequence of $\ell_\vv \circ G \circ f_\nu$ has itself a subsequence
that converges uniformly on the compact set $Y$. Consequently we
have
$$
\langle \yy^* , \vv \rangle = 0, ~~~\forall \vv \in \cH.
$$
However this is impossible, since $\yy^* \in \Sigma$ and since
$\Sigma$ does not contain the zero vector.
\qed
\bigskip  \\

Next, for the strongly pseudoconvex boundary point $\pp$ and the
automorphisms $f_j$ of $\O$ as in ($\dag$), we show:

\begin{lemma} \sl
Let $\xx_1, \xx_2 \in \O$ be chosen such that there exists $r>0$
satisfying
$$
t \xx_1 + (1-t) \xx_2 \in \O, \qquad \forall t \hbox{ with } -r < t <
1+r.
$$
If $\displaystyle{\lim_{j\to\infty} f_j (\xx_1) = \pp}$, then it
follows that $\displaystyle{\lim_{j\to\infty} f_j (\xx_2) =
\pp}$.
\end{lemma}

\noindent
\bf Proof: \rm
First of all, because of the strong pseudoconvexity, we may choose
an open neighborhood $V$ with $\pp \in V \subset U$ on which there
exists a continuous \it local peak function \rm $g : V \cap
\overline{\O} \to \CC$ satisfying
\medskip \\

\begin{itemize}
\item[(A)] $g$ is continuous on $V \cap \overline{\Omega}$,
\item[(B)] $g$ is holomorphic on $V \cap \O$,
\item[(C)] $g(\pp) = 1$, and
\item[(D)] $|g(\ww)| < 1$ for every
   $\ww \in V \cap \overline{\O} \setminus \{\pp\}$.
\end{itemize}
\vspace*{.15in}
Indeed, one directly imitates the construction and properties of the
Levi polynomial $L_P(z)$ (see [KRA1, p.\ 212]); then $\hbox{exp}(L_P(z))$
does the job.
\medskip  \\

We will shrink $U$ so that $U=V$ from here on.  In particular, $g$
is now defined on $U$.
\medskip \\

Let us write $D_s = \{ z \in \CC \mid |z| < s \}$ for each
positive real number $s$, and let $D = D_1$ denote the open unit disc.
Complexifying the straight line segment through $\xx_1$
and $\xx_2$, and then applying the Riemann mapping theorem if
necessary, we have an injective holomorphic map $h :
D_{(1+\delta)} \to \Omega$ such that $h(0)= \xx_1$ and $h(\lambda)
= \xx_2$, for some $\lambda \in D$ and $\delta > 0$.
\medskip  \\

The preceding lemma gives us an open neighborhood
$W_{\xx_1}$ of $\xx_1$ in $\O$ such that there exists $N_1 > 0$
satisfying the condition $f_j (W_{\xx_1}) \subset U \cap \O$
for all $j > N_1$. Note that there exists $\epsilon>0$ such that
$$
h(D_\epsilon) \subset W_{\xx_1}.
$$
Now consider the sequence of mappings
$$
g \circ f_j \circ h |_{D_\epsilon} : D_\epsilon \to D
$$
for $j =1,2, \ldots$.  Notice that every subsequence has
itself a subsequence that converges uniformly on compact subsets to a
holomorphic mapping that has value $1$ at $0$.  Thus
(by the maximum principle) the limit
mapping is the constant mapping $1$, and hence the above sequence
itself converges to the constant mapping $1$, uniformly on compact
subsets.  In particular, we obtain that
$$
\lim_{j \to \infty} \sup_{z\in D_r} \|f_j \circ h (z) - \pp \| = 0
$$
for every $0 < r< \epsilon$.
\medskip \\

At this point, we will use the quadratic polynomial mapping
$G:\cH \to \cH$ and the hyperplane $\Sigma$ used in the proof of
the preceding lemma.  Recall that $G$ maps $U \cap \O$ to a
strongly convex set and $\Sigma$ separates $G(U \cap \O)$ into two
disjoint subsets.
\medskip  \\

If $f_j (\xx_2)$ does not converge to $\pp$ in norm (as $j \to
\infty$), then we may move $\Sigma$ closer to $\pp$ if necessary
and then choose a subsequence of $f_j (\xx_2)$ (which,
by an abuse of notation, we denote in the
same way), so
that $f_j (\xx_2)$, for every $j$, belongs to the connected
component of $G(U \cap \O) \setminus \Sigma$ which does not contain
$\pp$.  Therefore, considering the image of the straight line
segment joining $0$ and $\lambda$ in $D$, one obtains points $\xi_j \in D$
such that $G \circ f_j\circ h(\xi_j) \in \Sigma \, \cap \, \O$ for every
$j$.  Applying the Banach-Alaoglu theorem again, we may assume
(choosing a subsequence if necessary) that $G \circ f_j \circ h
(\xi_j)$ converges weakly to $\yy^*$.  Again, the convexity
implies that $\yy^* \in \overline{\Sigma \cap \O}$.
\medskip  \\

Let $\vv \in \cH$, and write $\ell_\vv (\zz) = \langle \zz,
\vv \rangle$ for every $\zz \in \cH$. Consider the mappings
$$
\ell_v \circ G \circ f_j \circ h : D_{(1+\delta)} \to \CC
$$
for $j=1,2, \ldots$.  Notice that all their images are contained
in the bounded domain $\ell_v \circ G (\O)$.  Therefore Montel's
theorem implies that a subsequence can be chosen (for which we use the
same notation again) which converges uniformly on the
closed disc $\overline{D}$.  But we already know that $G \circ
f_j \circ h$ converges to zero uniformly on $D_{\epsilon/2}$.
Therefore we see that $\ell_v \circ G \circ f_j \circ h$
converges uniformly on $\overline{D}$ to the zero function.  In
particular,
$$
\langle \yy^*, \vv \rangle = \lim_{j\to\infty} \langle G \circ f_j
\circ h (\xi_j), \vv \rangle = 0.
$$
Since $\vv \in \cH$ is arbitrary, we now have that $\yy^* =
\origin$.  However, this is absurd because $\Sigma$ does not
contain the origin of $\cH$.

This contradiction yields the assertion.
\qed
\smallskip \\

Since $\O$ is a connected open subset of $\cH$, it follows that
every point of
$\O$ can be joined to $\qq$ in ($\dag$) by a polygonal line in
$\O$ with finitely many break points.  Therefore we
immediately obtain the following

\begin{corollary} \sl
Assuming ($\dag$) for $\O$, it follows that
$$
\lim_{j\to\infty} f_j (\zz) = \pp
$$
for every $\zz \in \O$.
\end{corollary}

Now we present the following strengthened form of the localization
lemma.

\begin{lemma} \sl
Let $\Omega$ be a bounded domain in $\cH$ with a boundary point
$\pp \in \bo$ such that
\begin{itemize}
\item[(1)] $\bo$ is $C^2$ smooth and strongly pseudoconvex near $\pp$, and
\item[(2)] there exist $\qq \in \O$ and $f_\ell \in \auto$ such that
$f_\ell (\qq) \to \pp$ as $\ell \to \infty$ in norm.
\end{itemize}
Then there exists an increasing sequence of open subsets $\O_1
\subset \O_2 \subset \cdots$ of $\O$ satisfying:
\begin{itemize}
\item[(\romannumeral1)]
$\displaystyle{\bigcup_{j=1}^\infty \O_j = \Omega}$,
\item[(\romannumeral2)]
for every open neighborhood $U$ of $\pp$ in $\cH$, there exists a
subsequence $f_{\ell_j}$ of $f_\ell$ above such that $f_{\ell_j}
(\O_k) \subset \O \cap U$ for every $j \ge k$.
\end{itemize}
\end{lemma}

\bf Proof: \rm We consider a countable dense subset (in the norm
topology)
$$
Z = \{\xx_j \in \O \mid j=1,2,\ldots\}
$$
of $\Omega$.  For each $j$, the preceding lemmata allow us to choose,
using the usual diagonalization process of selecting subsequences,
a sequence of real numbers $r_j \equiv r_{x_j} > 0$ and positive
integers $\nu_j$ such that
$$
f_\nu (\nbd{\xx_j}{r_j}) \subset U \cap \O, ~~\forall \nu > \nu_j.
$$
Thus we set
$$
\ell_j \equiv \max\{\nu_1, \dots, \nu_j\} + j
\hbox{ for each } j=1,2,\ldots
$$
so that
$$
\ell_1 < \ell_2 < \cdots
$$
and
$$
f_{\ell_j} (\nbd{\xx_k}{r_k}) \subset U \qquad \forall j \ge k.
$$
Now it suffices to let $$ \O_k \equiv \bigcup_{\ell=1}^k
\nbd{\xx_\ell}{r_\ell} $$ for each $k=1,2,\ldots$.  This concludes
the proof of the lemma as well as the construction
of the localization process.
\qed

\section{The Scaling Process}

We now discuss the scaling method (originated by S. Pinchuk in
$\CC^n$---see [IK1] for the history) and how it should
be adjusted so that it is a useful
technique in the study of automorphisms of bounded domains in
infinitely many dimensions.

For this purpose, we continue to work with an
orthonormal basis $\ee_1,\ee_2,\ldots$
for the separable Hilbert space $\cH$ under consideration.

\subsection{Normalization at the strongly pseudoconvex point}

We continue our exposition with the notation above; let $\O$ be
a bounded domain in $\cH$ with a $C^2$ strongly pseudoconvex
boundary point $\pp$ at which the automorphism orbit $f_j (\qq)$
accumulates.
\medskip  \\

{\it Furthermore, throughout the rest of the paper we utilize the
following notation: for each $\zz = \displaystyle{\sum_{j=1}^\infty z_k
\ee_k} \in \cH$, we set}
$$
\zz' = \sum_{j=2}^\infty z_k \ee_k.
$$
\medskip  \\

As before, let $U$ be an open neighborhood of $\pp$ and let
$G:U \to \cH$ be the quadratic change of coordinates which
maps $U$ biholomorphically onto $G(U)$ in such a way that
\begin{itemize}
\item $G({\bf p})= \origin$, where $\origin$ denotes the
origin of $\cH$, and
\item $\O_U \equiv G(\O \cap U)$ is represented by
$$
\O_U = \left \{ \zz = \sum_{j=1}^\infty z_k \ee_k \in G(U)
 \, \big| \,
\Re z_1 > \psi \left( (\Im z_1) \ee_1 + \zz' \right) \right \}
$$
\end{itemize}
where $\psi : \cH \to \RR$ is a $C^2$ function satisfying

\begin{itemize}
\item[{\boldmath $\circ$}]
$\psi (\origin) = 0$;
\item[{\boldmath $\circ$}]
$d \psi (\origin;\vv) = 0$ for every $\vv \in \cH$;
\item[{\boldmath $\circ$}]
there exists a constant $C_0 > 0$ such that
$$
d^2 \psi (\xx;\vv,\vv) \ge C_0 \| \vv \|^2, ~\forall \xx \in G(U), \
~\forall \vv \in \cH.
$$
\end{itemize}

We write
$$
\rho_U (\zz) = - \Re z_1 + \psi \left( (\Im z_1) \ee_1 + \zz' \right) ,
$$
which is the defining function of $\O_U$ in $G(U)$; we are going
to use $\rho_U$ repeatedly in the arguments below.

\subsection{Centering of the automorphism orbit}

Let us fix the notation for the automorphism orbit: put
$$
\qq_\nu \equiv G \circ f_\nu (\qq) \hbox{ for } \nu = 1,2, \ldots.
$$
(We may need to pass to a subsequence of $f_\nu$ when appropriate.)
\medskip \\

Consider the sequence of points $\qq_\nu \in \O_U$ approaching
the origin $\origin$ as $\nu \to \infty$.  First of all, we write
$$
\qq_\nu = \sum_{j=1}^\infty q_{\nu j} \ee_j, ~~q_{\nu j} \in \CC.
$$
Then, for each $\nu=1,2,\ldots$, we consider the boundary point
$$
\pp_\nu \equiv q_{\nu 1}^* \ee_1 + \qq'_\nu.
$$
Here $q_{\nu 1}^* \in \CC$ is chosen to satisfy the relations
$$
\pp_\nu \in \bo_U \hbox{ and } q_{\nu 1} - q_{\nu 1}^* > 0
$$
for each $\nu = 1,2,\ldots$.

Then we apply the complex affine linear isomorphism
$\Psi_\nu:\zz \mapsto \ww : \cH \to \cH$ defined by
\begin{eqnarray*}
w_1 & = &
  e^{i\theta_\nu} (z_1 - q_{\nu 1}^*) + T_\nu (\zz' - \qq'_\nu) \\
\ww' & = & \zz' - \qq'_\nu ,
\end{eqnarray*}
where the number $\theta_\nu \in \RR$ and the linear functional
$T_\nu : \cH \to \CC$ are to be chosen to satisfy the following conditions:

\begin{itemize}
\item The transformed domain $\Psi_\nu (\O_U)$ is supported by the hyperplane
defined by $\Re w_1 = 0$, and
\item $\Psi_\nu (\O_U) \subset \left \{\ww = \sum_j w_j {\bf e}_j
\mid \Re w_1 > 0 \right \}$
\end{itemize}
These conditions are easily satisfied if we choose $\theta_\nu$ and
$T_\nu$ as specified below.  Recall that
$\rho_U = - \Re w_1 + \psi (\Im w_1, \ww')$
is the defining function for $\O_U$.  Our choices for
$\theta_\nu$ and $T_\nu$ are determined by the following equations:
\begin{eqnarray*}
\theta_\nu
  & = & \tan^{-1} \left( \frac{\partial \psi}{\partial(\Re w_1)}
    \right) \ , \\
T_\nu \left(\sum_{k=2}^\infty v_k \ee_k \right)
  & = & 2e^{i\theta_\nu}
  \sum_{k=2}^\infty \frac{\partial\psi}{\partial w_k}
    (\Im q_{\nu 1}^*,\qq'_\nu) \, v_k \ ,
\end{eqnarray*}
and
$$
T_\nu(\ee_1) = 0 \, .
$$
Notice that $T_\nu$ is a bounded linear functional on $\cH$.
\medskip  \\

In conclusion, we have achieved that $T_\nu (\O_U)$ satisfies the
following properties (for every $\nu=1,2,\ldots$):

\begin{itemize}
\item[(C1)]
  $\Psi_\nu (\O_U)$ is supported by the real hyperplane defined
  by $\Re w_1 = 0$ at the origin $\origin$ in $\cH$;
\item[(C2)]
  $\origin \in \partial (\Psi_\nu (\O_U))$;
\item[(C3)]
  $\Psi_\nu (\qq_\nu) = e^{i\theta_\nu} (q_{\nu 1} - q_{\nu_1}^*) \ee_1$.
\end{itemize}
Also, notice that $e^{i\theta_\nu}$ converges to $1$ as $\nu$ tends
to $\infty$.

\subsection{The scaling map for $\O$ along $f_\nu$}

Now we consider the linear isomorphism $L_\nu : \cH \to \cH$
defined by
$$
L_\nu (\ww) = \frac{w_1 \ee_1}{\lambda_\nu}
+ \frac{\ww'}{\sqrt{\lambda_\nu}}
$$
for each $\nu$, where
$$
\lambda_\nu = q_{\nu 1} - q_{\nu 1}^*.
$$
The \it scaling process along $\qq_\nu$ \rm consists of the sequence
of maps defined by
$$
\sigma_\nu \equiv L_\nu \circ \Psi_\nu \circ G \circ f_\nu :
  \O_\nu \to \cH,  \hbox{ for } \nu = 1,2, \ldots
$$
In what follows, we shall develop properties of the scaling process,
and we also shall implement several refinements.

\section{Initial Scaling at a Strongly Pseudoconvex Point}

The goal of this section is to show that the scaling sequence
$\{\sigma_\nu \mid \nu=1,2,\ldots \}$ admits a subsequence that
converges to a mapping that maps $\O$ biholomorphically onto
a domain of $\cH$ which is in turn biholomorphic to the unit
ball in $\cH$.

\subsection{Ellipsoidal envelopes}

Let us return to $\O_U$ with the defining function
$\rho_U (\zz) = - \Re z_1 + \psi ( (\Im z_1)\ee_1 + \zz')$
which is $C^2$ smooth at every point of $G(U)$; the defining function
admits a constant $C_0 > 0$ satisfying the condition
$$
d^2\rho_U (\origin; \vv, \vv) \ge  C_0 \|\vv\|^2 \qquad \forall \vv \in \cH.
$$

Recall the centering map $\Psi_\nu$ from Section 5.2
associated with the point
sequence $\qq_\nu = G \circ f_\nu (\qq)$.  Let us put
\begin{eqnarray*}
\rho_\nu (\zz) & \equiv & \rho \circ \Psi_\nu^{-1} (\zz) \\
\Phi_\nu^{-1} (\zz) & \equiv & \Psi_\nu^{-1} (\zz) - \pp_\nu
\end{eqnarray*}
Notice that $d\Psi_\nu^{-1} (\xx; \zz) = \Phi_\nu^{-1} (\zz)$ for all
$\zz \in \cH$, regardless of the location of $\xx$ in $\cH$.  See
the previously defined
notation (Subsection 5.2) for $\pp_\nu$ as well as the definition
of the centering map $\Psi_\nu$.  Now
a direct computation yields
\begin{eqnarray*}
d\rho_\nu (\qq;\vv) & = & d\rho (\Psi_\nu^{-1} (\qq); \Phi_\nu^{-1} (\vv) ) \\
d^2\rho_\nu (\origin; \vv, \vv) \, ,
  & = & d^2 \rho (\pp_\nu; \Phi_\nu^{-1} (\vv), \Phi_\nu^{-1} (\vv) ) \, .
\end{eqnarray*}
Then we see that
\begin{eqnarray*}
& & \frac{1}{\|\zz\|^2} \big[ \rho_\nu (\zz) - \rho_\nu (\origin)
  - d\rho_\nu (\origin ; \zz) \big] \\
&=& \frac{\rho_U (\Psi_\nu^{-1}(\zz)) - \rho_U (\Psi_\nu^{-1} (\origin))
   - d\rho_U (\Psi_\nu^{-1}(\origin);\Phi_\nu^{-1} (\zz))}
     {\|\Phi_\nu^{-1} (\zz)\|^2} \cdot
    \frac{\|\Phi_\nu^{-1} (\zz)\|^2}{\|\zz\|^2} \\
\end{eqnarray*}

Recall that the sequence $\Psi_\nu^{-1}$ converges to the identity map of $\cH$
uniformly on every bounded subset of ${\mathcal H}$
and that $\rho_U$ is $C^2$ smooth. Notice
that we may choose $N>0$ and $r>0$ such that
$$
\|\Phi_\nu^{-1} (\zz) \| \ge \frac12 \|\zz\|,  \qquad \forall \zz \in \cH
\hbox{ with } \|\zz\| \le r \hbox{ and } ~\forall \nu > N.
$$
Then, shrinking $r>0$ if necessary, we see that the $C^2$ smoothness
of $\rho_U$ implies that
$$
\frac{\rho_U (\Psi_\nu^{-1}(\zz)) - \rho_U (\Psi_\nu^{-1} (\origin))
   - d\rho_U (\Psi_\nu^{-1}(\origin);\Phi_\nu^{-1} (\zz))}
      {\|\Phi_\nu^{-1} (\zz)\|^2} \cdot
    \frac{\|\Phi_\nu^{-1} (\zz)\|^2}{\|\zz\|^2}
$$
converges to
$$
d^2 \rho_U \left (\origin ; \frac{\zz}{\|\zz\|}, \frac{\zz}{\|\zz\|}\right )
$$
uniformly on $\nbd{\origin}{r}$ in $\cH$ as $\nu$ tends to infinity.
Therefore we are able to deduce the following:

\begin{lemma}  \sl
There exist positive constants $r$, $N$, and $c$
satisfying the inequality
$$
\rho_\nu (\zz) - \rho_\nu (\origin) - d\rho_\nu (\origin ; \zz)
> c\|\zz\|^2
$$
for every $\zz \in \cH$ with $\|\zz\|<r$ and for every $\nu > N$.
\end{lemma}

\subsection{First adjustment of the scaling sequence}

Recall that the domain $L_\nu \circ \Psi_\nu (\Omega_U)$
is defined by the inequality
$$
c_\nu \Re z_1 > d^2 \rho_\nu (\origin; \zz', \zz') + R_\nu (\zz) \, ,
$$
where:

\begin{itemize}
\item $c_\nu$ converges to $1$;
\item $R_\nu (\zz)$ tends to zero uniformly on each bounded set
\end{itemize}
as $\nu$ tends to $\infty$.

For each sufficiently large $\nu$, we consider the linear
transformation $S_\nu : \cH \to \cH$ satisfying

\begin{itemize}
\item $S_\nu (\ee_1) = \ee_1/c_\nu$,
\item $S_\nu \left( (\ee_1)^\perp \right) \subseteq (\ee_1)^\perp$, and
\item $d^2 \rho_\nu (\origin; S_\nu^{-1} (\zz'), S_\nu^{-1} (\zz'))
= \langle \zz' , \zz' \rangle_\cH$.
\end{itemize}
Notice that we can choose $S_\nu$ so that both $S_\nu$ and $S_\nu^{-1}$
are uniformly bounded linear operators on $\cH$.

Notice that there exists $r>0$ such that
$$
\langle S_\nu^{-1}(\zz'), S_\nu^{-1}(\zz') \rangle \ge
r \langle \zz', \zz' \rangle , \qquad \forall \zz' \in (\ee_1)^\perp \, ,
$$
for every $\nu = 1,2, \ldots$.
\medskip  \\

Therefore we have the following:

\begin{quote}
\begin{itemize}
\item[\bf (i)] $S_\nu \circ L_\nu \circ \Psi_\nu (\Omega_\nu)$ is defined
by
$$
\Re z_1 > \langle \zz', \zz' \rangle + R_\nu (\zz) \, ,
$$
where $R_\nu$ converges to zero uniformly on every bounded
subset of $\cH$, as $\nu \to \infty$.

\item[\bf (ii)] $S_\nu \circ L_\nu \circ \Psi_\nu (\O_U) \subset
\cE (1,r)$ for some fixed $r>0$ independent of $\nu$, where
$$
\cE (1,r)
= \{ \zz \in \cH \mid \Re z_1 > r \langle \zz', \zz' \rangle \}.
$$
\end{itemize}
\end{quote}
\vspace*{.15in}

Then we consider the linear fractional transformation
$F : \cH \to \cH$ defined by
$$
F(\zz) = \left( \frac{z_1 - 1}{z_1 + 1} \right) \ee_1
+ \frac{2}{z_1 + 1} \sum_{k=2}^\infty z_k \ee_k.
$$
This mapping transforms $\cE (1,r)$ biholomorphically onto
$\BB_r \equiv \{ \zz \in \cH \mid |z_1|^2 + r \| \zz'\|^2 < 1 \}$.
\bigskip  \\

Finally, we write
\begin{eqnarray*}
\omega_\nu & \equiv & F \circ S_\nu \circ \sigma_\nu  \\
  & = & F \circ S_\nu \circ L_\nu \circ \Psi_\nu \circ G \circ f_\nu
\end{eqnarray*}
for each $\nu = 1,2, \ldots$.  The map $\omega_\nu$ sends $\Omega_\nu$ into
$\BB_r$ for every $r > 1$.

\subsection{Convergence of the First Adjustment}

We first specify the domains $\O_\nu$ which exhaust $\Omega$.
\medskip  \\

Let us consider
$$
\Sigma_n = \CC \ee_1 \oplus \cdots \oplus \CC \ee_n
~~~n = 1,2, \ldots
$$
Let $Z_k$ be a countable dense subset of $\Sigma_k$ for each
$k =1,2,\ldots$, and let
$$
Z = \bigcup_{k=1}^\infty Z_k.
$$
Let us write $Z = \{\xx_1, \xx_2, \ldots \}$.
\medskip \\

Then, as in the localization process introduced in Section 4,
we may choose the balls $\nbd{\xx_j}{r_j}$ and define
$$
\O_\ell = \bigcup_{j=1}^\ell \nbd{\xx_j}{r_j}
$$
for each $\ell = 1,2, \ldots$.

Let $r>0$ be chosen arbitrarily.  Then notice that the steps {\bf (i)}
and {\bf (ii)} of Subsection 6.2
will be valid for $\cE (1,r)$ if we choose $U$ sufficiently small.
Therefore
the localization lemma in Section 4 allows us to choose a subsequence so that
$$
\omega_\nu (\Omega_\nu) \subset \BB_r
$$
for every $\nu = 1,2, \ldots$.
\medskip  \\

Now the Banach-Alaoglu theorem allows us to choose a subsequence
$\omega_{\nu_k}$ of $\omega_\nu$ such that $\omega_{\nu_k} (\xx_1)$
converges weakly to $\widehat\omega (\xx_1)$ as $\nu \to \infty$.
Then we can inductively choose a subsequence $\{\omega_{\nu_k}\}$
of $\omega_{\nu_{k-1}}$ ~so that
$\omega_{\nu_k} (\xx_j)$ converges weakly to $\widehat\omega (\xx_j)$
as $\nu \to \infty$ for every $j = 1,2,\ldots, k$.

Re-naming the sequence $\omega_\nu = \omega_{\nu_\nu}$, we
see that
$$
\omega_\nu (\xx) \to \widehat\omega (\xx)
\hbox{ weakly, as } \nu \to \infty \, ,
$$
for every $\xx \in Z$.
\medskip \\

Let $\ell \in \cH^*$. Consider the sequence
$$
\ell \circ \omega_\nu |_{\Sigma_n \cap \O_\nu}
: \Sigma_n \cap \O_\nu \to \CC \, .
$$
Notice that this is a uniformly bounded sequence. By Montel's theorem,
there exists a subsequence $\ell \circ \omega_{\nu_k}$ that converges
to the holomorphic function, say $\phi: \Sigma_n \cap \O \to \CC $,
uniformly on compact subsets of $\Sigma_n \cap \Omega$. Notice that
for each $\xx \in Z_n$ we have
\begin{eqnarray*}
& & \hspace{-35pt}
|\phi (\xx) - \ell \circ \widehat\omega (\xx)|  \\
& \le &
|\phi(\xx) - \ell \circ \omega_{\nu_k}(\xx)|
+ |\ell \circ \omega_{\nu_k} (\xx) - \ell \circ \widehat\omega (\xx)|  \\
& \longrightarrow & 0
\end{eqnarray*}
as $k \to \infty$.  Therefore $\ell \circ \widehat\omega = \phi$
on $Z_n$.  In particular, $\ell \circ \widehat\omega$ extends uniquely to
a holomorphic function on $\Sigma_n \cap \Omega$, which we will
denote with the same notation.
\medskip  \\

Now we prove

\begin{lemma} \sl
The sequence $\omega_\nu (\zz)$ converges weakly
to an element $\widehat{\omega}(\zz)$ for every
$\zz \in \Sigma_n \cap \O$.
\end{lemma}
{\bf Proof:} \rm
Recall that $\omega_\nu$ forms an equibounded family of holomorphic
mappings.  Therefore the Banach-Alaoglu theorem allows us to choose
a subsequence $\omega_{\nu_k}$ such that
$$
\omega_{\nu_k} (\zz) \to \alpha_\zz \in \cH \hbox{ weakly}
$$
as $k \to \infty$.

Let $\ell \in \cH^*$.  Applying Montel's theorem to the restriction
of these maps to $\Sigma_n \cap \O$, we may again choose a
subsequence $\omega_{\nu_{j_k}}$ such that
$$
\ell \circ \omega_{\nu_{j_k}} (\zz) \to \psi(\zz)
$$
uniformly on compact subsets of $\Sigma_n \cap \O$.  Notice
that we then obtain
$$
\psi (\xx)
= \ell \circ \widehat\omega (\xx) ~~~\hbox{ for every } \xx \in Z_n
$$
by the argument preceding the statement of this lemma.  Therefore
we can easily deduce that $\psi (\zz)$ is indeed uniquely determined to
be  $\ell \circ \widehat\omega(\zz)$, for every $\ell \in \cH^*$.
Consequently we have
$$
\alpha_\zz = \sum_{k=1}^\infty \ell_k \circ \widehat\omega (\zz) \ee_k,
$$
where $\ell_k (\ww) = \langle \ww, \ee_k \rangle$ for every $\ww \in \cH$.

Hence it suffices to let $\widehat\omega(\zz) = \alpha_\zz$
for each $\zz \in \Sigma_n \cap \Omega$.
\qed
\bigskip \\

Note that the lemma defines $\widehat\omega$ on $\Sigma_n \cap \Omega$.  Since
$\ell \circ \widehat\omega$ is a holomorphic function on $\Sigma_n \cap \Omega$
for every $\ell \in \cH^*$, Hartogs's analyticity theorem implies:

\begin{lemma}  \sl
The restriction
$\widehat\omega|_{\Sigma_n \cap \Omega}: \Sigma_n \cap \Omega \to \cH$
is a holomorphic mapping for every finite-dimensional slice
$\Sigma_n \cap \Omega$ of $\Omega$.
\end{lemma}

Now we extend $\widehat\omega$ to a holomorphic mapping of $\Omega$
into $\cH$.  Let
$$
\zz = \sum_{k=1}^\infty z_k \ee_k \in \O.
$$
We write
$$
\zz_n = \sum_{k=1}^n z_k \ee_k
$$
for every $n =1,2, \ldots$.
Note that $\zz_n \in \Sigma_n \cap \Omega$, and $\zz_n \to \zz$
in norm as $n \to \infty$.
\medskip \\

Let $Q \equiv \{\zz\} \cup \{\zz_n \mid n=1,2, \ldots\}$.
Let $\ell \in \cH^*$ be arbitrarily given.  Since
$Q$ is compact, it follows that a subsequence $\omega_{\nu_j}$ of $\omega_\nu$
can be chosen so that $\ell \circ \omega_{\nu_j}$ converges
to $\ell \circ \widehat\omega$, uniformly on $Q$.  Therefore
it follows immediately that
the sequence $\ell \circ \widehat\omega (\zz_n)$ is uniformly Cauchy.
Thus it converges---since $\cH$ is complete.  Let us write
$$
\beta_k \equiv
\lim_{n\to\infty} \langle \widehat\omega (\zz_n), \ee_k \rangle
$$
for each $k = 1,2, \ldots$.  Then we define
$$
\widehat\omega (\zz) = \sum_{k=1}^\infty \beta_k \ee_k.
$$
This allows us to extend $\widehat\omega$ to all of $\Omega$.
\medskip \\

Now the question we need to resolve is whether the mapping $\widehat\omega$
so defined is a holomorphic mapping of $\Omega$ into $\cH$.
To answer this question we have:

\begin{proposition} \sl
The mapping $\widehat\omega : \Omega \to \cH$ is holomorphic, and
$\widehat\omega (\Omega) \subset \BB$.
\end{proposition}

\noindent \bf Proof. \rm Since
\begin{multline*}
|\ell \circ \omega_\nu (\zz ) - \ell \circ \widehat\omega (\zz)|
\le
|\ell \circ \omega_\nu (\zz) - \ell \circ \omega_\nu (\zz_n)| \\
+
|\ell \circ \omega_\nu (\zz_n) - \ell \circ \widehat\omega (\zz_n)|
+
|\ell \circ \widehat\omega (\zz_n) - \ell \circ \widehat\omega (\zz)| \, ,
\end{multline*}
we can deduce immediately that $\omega_\nu$ converges weakly to
$\widehat\omega$ at every point of $\O$.

Now consider the $\cH$-valued mapping
$$
g (\lambda) \equiv \widehat\omega (\pp+\lambda \vv)
$$
of $\lambda \in \CC$, where $\pp \in \O$ and $\vv \in \cH$.  The function
$g$ is well-defined for all $\lambda \in \CC$ with $|\lambda|<r$ for
some sufficiently small $r>0$.

Let us consider the projections $\pi_n : \cH \to \CC\ee_1 \oplus
\ldots \oplus \CC\ee_n$ for each $n =1,2, \ldots$ and examine
the mapping
$$
g_n (\lambda) \equiv \widehat\omega (\pi_n(\pp) + \lambda \pi_n (\vv)) \, ,
$$
which is well-defined for all $\lambda \in \CC$ with
$|\lambda|< r$. (Shrink $r>0$ if necessary.  But note that one may choose
$r$ independent of $n$.)

Lemma 6.3 implies that $g_n$ is holomorphic for every $n$, and
$\{g_n\}$ forms an equibounded family.  Therefore, for each $\ell \in \cH^*$,
$\ell \circ g_n$ admits a subsequence, say $\ell \circ g_{n_j}$,
that converges to a holomorphic function, say $h$, uniformly
on compact subsets.  On the other hand, the preceding arguments
yield that
$\ell \circ g_n (\lambda)
= \ell \circ \widehat\omega (\pi_n (\pp + \lambda\vv))$ converges to
$\ell \circ \widehat\omega (\pp + \lambda\vv)$ pointwise, by
construction.  Therefore we see that
$$
\ell \circ \widehat\omega (\pp + \lambda\vv) = h(\lambda)
$$
for every $\lambda \in \CC$ with $|\lambda| < r$. Since $\ell
\in \cH^*$ is arbitrary, we obtain that $\widehat\omega (\pp + \lambda\vv)$
is a holomorphic mapping of $\lambda$ for every fixed $\pp \in \Omega$
and $\vv \in \cH$.

Notice that the weak limit of a sequence in $\cH$ stays in the closed
convex hull of the given sequence.  Therefore we note that
the set $\widehat\omega (\O)$ is bounded.  Altogether, the mapping
$\widehat\omega : \O \to \cH$ is bounded, and holomorphic on every
one-dimensional slice (and hence in every finite-dimensional slice
by Hartogs's analyticity theorem).  Consequently, we may conclude
that $\widehat\omega: \O \to \cH$ is holomorphic. (Notice that
the Gateaux definition of holomorphicity formally requires $C^1$ smoothness,
but the local boundedness and analyticity in every
finite-dimensional slice give sufficient grounds for holomorphicity
of the function.  See for instance [Mu].)

Now we repeat that one may adjust the subsequence $\omega_\nu$
(choosing a subsequence if necessary) so that
$\omega_\nu (\O_\nu) \subset (1+\frac1\nu ) \BB$.  We have addressed
this point in Section 6.3 above, in which the localization lemma
was essential.

Using again the fact that the weak limits stay in the closed convex
hull of the sequence in consideration, we obtain finally that
$$
\widehat\omega (\O) \subset \overline{\BB}.
$$
Recall, by definition, that $\omega_\nu (\qq) = \origin$ for every
$\nu = 1,2, \ldots$.  Thus $\widehat\omega (\qq) = \origin$.  Then
we may apply the Maximum Principle as in the finite-dimensional case
to conclude immediately that $\widehat\omega (\O) \subset \BB$.
This completes the proof.
\qed
\medskip \\

It is also true that the same principle as in the proof of Lemma 6.4
yields a holomorphic weak limit mapping of the inverse maps of the
elements of the above scaling sequence.  Just the same, it is not yet
enough to conclude that the scaled limit mapping so constructed above
is a biholomorphism.  So we shall adjust the scaling process again in
the next section.

\section{Second Adjustment of Scaling}

We begin this section with some functional-analytic arguments.

\subsection{Base-change with Uniformly Bounded Operators}

Consider a sequence $T_\nu : \cH \to \cH$ of invertible
$\CC$-linear operators admitting a constant $C>0$ such that
$$
\frac1C \| \vv \| \le \| T_\nu (\vv) \| \le C \| \vv \|
\hbox{ and }
\frac1C \| \vv \| \le \| T_\nu^{-1} (\vv) \| \le C \| \vv \|
$$
for every $\vv \in \cH$ and $\nu=1,2, \ldots$.
In general, we only have the subsequential weak convergence for
a sequence of such operators.  The momentary goal is to improve
its convergence by composing with Hilbert space linear isometries.

\subsubsection{Gram-Schmidt Process}

Consider $\ee_1, \ee_2, \ldots$, a fixed orthonormal basis
system of $\cH$.  We perform the Gram-Schmidt process as
follows:
\begin{eqnarray*}
\ff_{\nu 1} & := & T_\nu (\ee_1)  \\
        & \vdots &  \\
\ff_{\nu k} & := &
T_\nu (\ee_k)
- \sum_{j=1}^{k-1} \frac{\langle T_\nu (\ee_k), \ff_{\nu j} \rangle}
             {\langle \ff_{\nu j}, \ff_{\nu j} \rangle}  \ff_{\nu j} \\
        & \vdots &
\end{eqnarray*}
inductively on $k = 1,2, \ldots$.

\subsubsection{Uniform Upper-bound for $\ff_{\nu k}$}

By construction, we have that $\langle \ff_{\nu j}, \ff_{\nu k} \rangle = 0$
whenever $j \not= k$. Consequently,
$$
\| T_\nu (\ee_k) \|^2
= \|\ff_{\nu k} \|^2
+ \sum_{j=1}^{k-1}
 \left| \langle T_\nu (\ee_k), \ff_{\nu j}/\|\ff_{\nu j}\| \rangle \right|^2.
$$
In particular, we obtain

\begin{equation}   \label{one}
\|\ff_{\nu k} \| \le C
\end{equation}
and
\begin{equation}   \label{two}
\sum_{j=1}^{k-1}
 \left| \langle T_\nu (\ee_k), \ff_{\nu j}/\|\ff_{\nu j}\| \rangle \right|^2
\le C
\end{equation}
for all positive integer values of $\nu$ and $k$.

\subsubsection{Uniform lower-bound for $\ff_{\nu k}$}

Note that our construction implies that
$$
\hbox{Span} \{\ff_{\nu 1}, \ldots, \ff_{\nu k} \}
 = \hbox{Span} \{T_\nu (\ee_1), \ldots, T_\nu (\ee_k)\}
$$
for every $k =1,2,\ldots$.  Hence
$$
\ff_{\nu k} = T_\nu (\ee_k)
+ \sum_{j=1}^{k-1} \lambda_{j} T_\nu (\ee_j)
$$
for some $\lambda_j \in \CC$, $j=1,2, \ldots$.  As a result,
\begin{eqnarray*}
\|\ff_{\nu k}\|^2
& = & \left\| T_\nu
  \left( \ee_k + \sum_{j=1}^{k-1} \lambda_{j} \ee_j ) \right) \right\|^2 \\
& \ge &
\frac1{C^2} \left \|\ee_k + \sum_{j=1}^{k-1} \lambda_{j} \ee_j  \right\|^2  \\
& = &
\frac1{C^2} \left( 1 + \sum_{j=1}^{k-1} |\lambda_{j}|^2 \right)
\end{eqnarray*}
\medskip \\

In conclusion, we have
\begin{equation} \label{three}
\|\ff_{\nu k} \| \ge \frac1C, \text{ for every } \nu, k = 1,2, \ldots.
\end{equation}

\subsection{A Sequence of Hilbert Space Isometries}

Consider the mapping $S_\nu : \cH \to \cH$ defined by
$$
S_\nu \left( \sum_{j=1}^\infty \alpha_j
\frac{\ff_{\nu k}}{\|\ff_{\nu k}\|} \right)
= \sum_{j=1}^\infty \alpha_j \ee_j
$$
for every sequence $\{\alpha_j\}_{j=1,2,\ldots}$ of complex numbers
satisfying
$$
\sum_{j=1}^\infty |\alpha_j|^2 < \infty.
$$
Notice that this defines a Hilbert space isometry
for every $\nu = 1,2, \ldots$.

\subsection{The sequence $S_\nu \circ T_\nu$}

Recall the notation
$$
\Sigma_n = \hbox{Span}\,\{\ee_1, \ldots, \ee_n\}, \hbox{ for }
n = 1,2, \ldots.
$$

Let us write
$$
A_\nu \equiv S_\nu \circ T_\nu
$$
for each $\nu$.  Then, as its construction shows, the sequence
$A_\nu$ satisfies:

\begin{itemize}
\item
$A_\nu$ maps $\Sigma_n$ onto itself as a $\CC$-linear isomorphism.
\item
Both $\|A_\nu|_{\Sigma_n} \|$ and $\|A^{-1}|_{\Sigma_n}\|$ are
bounded by the positive constant $C$, independent of $\nu$ and $n$.
\end{itemize}

Therefore we may choose a subsequence of $A_\nu$ which converges
in norm to an operator
$$
A : \bigcup_{n=1}^\infty \Sigma_n \to \bigcup_{n=1}^\infty \Sigma_n
$$
\it which is invertible on each $\Sigma_n$ \rm for $n=1,2,\ldots$.

We now extend $A$ to an operator on all of $\cH$.  Let $\epsilon > 0$ be
given.  For each $\vv \in \cH$ with $\|\vv\|=1$, we choose a
positive integer $k$ such that
$$
\|\vv - \pi_k (\vv) \| < \frac{\epsilon}{3C}.
$$
Here $\pi_k : \cH \to \Sigma_n$ denotes the orthogonal
projection, as before.  Notice, in particular,
that $\|\pi_k (\vv) \| \le \|\vv\| \le 1$.
Then, for this $k$, we may choose $N > 0$ such that
$$
\|A_\nu (\pi_k (\vv))  - A_\mu (\pi_k (\vv)) \|
\le \frac{\epsilon}{3}, \hbox{ for all } \nu, \mu > N.
$$
Then we obtain, for every $\nu, \mu > N$, that
\begin{eqnarray*}
\|A_\nu (\vv) - A_\mu (\vv)\| & \le & \|A_\nu (\vv-\pi_k (\vv)) \|
  +
\|A_\nu (\pi_k (\vv)) - A_\mu (\pi_k (\vv))\| \\
&& \qquad + \| A_\mu (\vv - \pi_k (\vv)) \| \\
& < & \epsilon.
\end{eqnarray*}
The completeness of $\cH$ shows that
$$
\widehat A (\vv) \equiv \lim_{\nu\to\infty} A_\nu (\vv)
$$
exists for every $\vv \in \cH$.  It is obvious that $\widehat A$
is linear and is the unique extension of $A$. Now we show that
$\widehat A$ is bounded.  To see this, let $\vv \in \cH$ satisfy
$\|\vv\| = 1$.  Then choose $\nu$ such that $\| \widehat A (\vv)
- A_\nu (\vv) \| \le 1$.  Thus
\begin{eqnarray*}
\| \widehat A (\vv) \| & \le &
 \|\widehat A (\vv) - A_\nu (\vv) \| + \| A_\nu (\vv) \| \\
& \le &
 1 + C.
\end{eqnarray*}

Repeating the same process with $A_\nu^{-1}$, and choosing subsequences
whenever necessary, we arrive at the following conclusion:

\begin{proposition} \sl
Every subsequence of $A_\nu$ above has itself a subsequence that converges
to a bounded operator $\widehat A : \cH \to \cH$, uniformly on each
finite slice $\Sigma_n = \hbox{\rm Span}\,\{\ee_1, \ldots, \ee_n\}$,
$n=1,2, \ldots$.  Moreover, $\widehat A$ is invertible and
preserves $\Sigma_n$ for each $n$.
\end{proposition}

Rephrasing, we have

\begin{proposition} \sl
Let $T_\nu : \cH \to \cH$ be invertible $\CC$-linear operators with
a constant $C > 0$ such that
$$
\|T_\nu\| \le C \quad \hbox{ and } \quad \|T_\nu^{-1}\| \le C
$$
for every $\nu = 1,2, \ldots$.  Then there exists a sequence
$S_\nu : \cH \to \cH$ of Hilbert space isometries such that
the sequence $A_\nu \equiv S_\nu \circ T_\nu$ satisfies the following:
\begin{itemize}
\item[(\romannumeral1)]
   Each $S_\nu \circ T_\nu$ maps $\Sigma_n$ isomorphically
   onto $\Sigma_n$.
\item[(\romannumeral2)]
  There exists a subsequence $A_{\nu_j}$
  such that both $A_{\nu_j}$ and $A_{\nu_j}^{-1}$ converge
  to invertible bounded operators, respectively, each of which maps
  $\Sigma_n$ isomorphically onto itself for every $n=1,2,\ldots$.
\end{itemize}
\end{proposition}

\subsection{The Second Adjustment of the Scaling Sequence}

We now apply the above arguments to the scaling sequence.  We begin with
the following observation on the derivatives of the scaling mappings.

\subsubsection{The Sequence $d\omega_\nu (\qq; \, \cdot \, )$}
We write $T_\nu( \, \cdot \, ) = d\omega_\nu (\qq; \cdot )$ for each $\nu$.  Then
we prove

\begin{proposition} \sl
There exists a subsequence (which we also denote by $T_\nu$) of
$T_\nu$ which admits a constant $C > 0$ such that
$$
\frac{\|\vv\|}{C} \le \|T_\nu (\vv) \| \le C \|\vv\|
$$
and
$$
\frac{\|\vv\|}{C} \le \|T_\nu^{-1} (\vv) \| \le C \|\vv\|
$$
for every $\vv \in \cH$.
\end{proposition}

\noindent \bf Proof. \rm
Recall that, choosing a subsequence, we may arrange that
\medskip  \\

\begin{itemize}
\item[(1)]
there exists a constant $\delta > 0$ independent of $\nu$ such that
$\BB (\qq; \delta) \subset \Omega_\nu$ for every $\nu$;
\item[(2)]
there exists a constant $r > 1$ independent of $\nu$ such that
$\omega_\nu (\O_\nu) \subset \BB(\origin; r)$;
\item[(3)]
$\omega_\nu (\qq) = \origin$.
\end{itemize}
\vfill
\eject

\noindent We use the notation $k_G : G \times \cH \to \RR$ to denote
the (infinitesimal) Kobayashi metric
of the domain $G \subset \cH$.  (Notice that the notion of Kobayashi metric
admits a natural generalization from the finite-dimensional case, at least
for our Hilbert space $\cH$.)  We immediately obtain
\begin{eqnarray*}
\frac1r \| T_\nu (\vv) \|_\cH
& = & k_{\BB(\origin; r)} (\origin; T_\nu (\vv))  \\
& = & k_{\BB(\origin; r)} (\omega_\nu (\qq); d\omega_\nu (\qq; \vv) ) \\
& \le & k_{\O_\nu} (\qq; \vv) \\
& \le & k_{\BB(\qq;\delta)} (\qq; \vv) \\
& = & \frac1{\delta} \|\vv\|.
\end{eqnarray*}
This yields that
$$
\| T_\nu \| \le \frac{r}{\delta}, \hbox{ for all } \nu = 1,2, \ldots.
$$
A similar argument on $T_\nu^{-1}$  yields a uniform bounds for
$\| T_\nu^{-1} \|$.  Altogether, we arrive at the desired conclusion.
\qed

\subsubsection{Adjusted Scaling Sequence $\alpha_\nu$}

Notice first that, for any subsequence of a scaling
sequence $\omega_\nu$, we may now choose a sequence
$S_\nu$ of $\CC$-linear isometries such that a subsequence
(again we abuse notation and let $\alpha_\nu$ also be
the subsequence)
$$
\alpha_\nu \equiv S_\nu \circ \omega_\nu : \Omega_\nu \to \cH
$$
satisfies the following conditions:

\begin{quote}
\begin{itemize}
\item[(\romannumeral1)]
For every $\nu = 1,2,\ldots$, the operator
$$
A_\nu \equiv d\alpha_\nu (\qq; \cdot):\cH \to \cH
$$
and its inverse $A_\nu^{-1}$
have their norms bounded by a constant independent of $\nu$;
\item[(\romannumeral2)]
Each $A_\nu$ maps $\Sigma_n = \hbox{Span}\, \{\ee_1, \ldots, \ee_n \}$
isomorphically onto itself, for every $n=1,2,\ldots$;
\item[(\romannumeral3)]
The sequence $A_\nu$ converges to a bounded invertible operator
(with a bounded inverse), say $\widehat A$, on each $\Sigma_n$
for $n=1,2,\ldots$.
\end{itemize}
\end{quote}

Notice that, for every $r>1$, there exists $N > 0$ such that
$\omega_\nu (\O_\nu) \subset \BB(\origin; r)$ for every $\nu > N$.
Since the $S_\nu$ are isometries, we obviously have the same property
for the newly adjusted scaling sequence $\alpha_\nu$.  Hence we
may choose a weak limit in the sense of Section 5, which eventually
yields the limit holomorphic mapping $\widehat\alpha:\O \to \BB$.
\medskip \\

Notice that
$$
d\widehat\alpha (\qq; \cdot ) = \widehat A.
$$
Choosing a subsequence once again if necessary, we see that there
exists a weak limit mapping $\beta: \BB \to \cH$ of the sequence
$$
\alpha_\nu^{-1} : \BB(\origin; s) \to \O
$$
for $s<1$.  The preceding arguments then show that
\begin{itemize}
\item[(\romannumeral1)]
$[d\alpha_\nu^{-1}](\origin; \cdot )$ maps each $\Sigma_n$ isomorphically
onto itself, and as a result it converges to ${\widehat A}^{-1}$ on
each finite-dimensional slice $\Sigma_n$.  Consequently, we have
$$
d\beta (\origin; \cdot) = {\widehat A}^{-1}.
$$
\item[(\romannumeral2)]
The limit mapping $\widehat\beta$ is well-defined on the unit ball $\BB$
and satisfies $\widehat\beta (\BB) \subset \widehat\O$, where $\widehat\O$
denotes the (open) convex hull of $\O$.
\end{itemize}
(In general, $\widehat\beta (\BB)$ is contained in the closure of
$\widehat\O$.  But, since $\widehat\beta (\origin) = \qq$, the Maximum
Principle implies the conclusion above.)
\medskip \\

Finally, we arrive at

\subsection{Proof of the Main Theorem}

We are now ready to give the proof of the Main Theorem (Theorem 3.1),
which is
\bigskip  \\

\noindent \bf Main Theorem. \it Every bounded convex domain $\O$ in a
separable Hilbert space $\cH$ admitting a $C^2$
strongly pseudoconvex boundary point at which an automorphism
orbit accumulates is biholomorphic to the open unit ball of $\cH$.
\bigskip \\

\noindent  \bf Proof. \rm
The scaling arguments above yields the holomorphic mappings
$$
\widehat\alpha : \O \to \BB
\hbox{ and }
\widehat\beta: \BB \to \O = \widehat\Omega
$$
satisfying the conditions
\begin{itemize}
\item[(1)]
$\widehat\alpha (\qq) = \origin$,
\item[(2)]
$\widehat\beta (\origin) = \qq$,
\item[(3)]
$d(\widehat\alpha \circ \widehat\beta) (\origin; \vv) = \vv$
and $d(\widehat\beta \circ \widehat\alpha) (\qq; \vv) = \vv$
for every $\vv \in \cH$.
\end{itemize}
Now we apply the infinite-dimensional version of Cartan's
Uniqueness Theorem (essentially identical with the original
finite-dimensional version, see [BM]).  Then both $\widehat\alpha \circ
\widehat\beta$ and $\widehat\beta \circ \widehat\alpha$ are
equal to the identity map.
Therefore $\alpha : \O \to \BB$ is in particular a biholomorphic
mapping.  This completes the proof.
\qed

\section{Concluding Remarks}

The main result of this paper---a version of the Wong/Rosay theorem
in a separable Hilbert space---has only been proved here for convex
domains.  The role of convexity in several complex variables has been
long established (see [LEM1], [KRA2]).  Nevertheless, it is not well
understood.  In particular, there is a rather poor understanding of
convexity from the biholomorphically invariant point of view.  We
have no examples to indicate whether our result ought to be true for
more general classes of domains.

Of course the exploration of domains of finite type will be an entirely
new world.  The situation for such domains in dimension 2 is now
fairly well in hand.  In higher dimensions there is much yet to
be understood, and we are still some distance from having any true grasp
of the infinite dimensional situation.

\vspace*{.35in}

Kang-Tae Kim

Department of Mathematics

Pohang University of Science and Technology

Pohang 790-784 The Republic of Korea

{\tt (kimkt@postech.edu)}
\bigskip  \\

Steven G. Krantz

Department of Mathematics

Washington University in St. Louis

St. Louis, MO 63130 U.S.A.

{\tt (sk@math.wustl.edu)}


\begin{thebibliography}{ABCDE}


\bibitem[BP1]{BP1}  Bedford, E.\ and Pinchuk, S.,
  Domains in ${\CC}^2$ with non-compact holomorphic
automorphism group (translated from Russian),
{\it  Math.\ USSR-Sb.} 63(1989), 141--151.\

\bibitem[BP2]{BP2}  Bedford, E.\ and Pinchuk, S.,
Domains in ${\CC}^{n+1}$ with
non-compact automorphism groups, {\it  J.\ Geom.\ Anal.} 1(1991),
165--191.

\bibitem[BP3]{BP3}  Bedford, E.\ and Pinchuk, S.,
Convex domains with non-compact automorphism group
(translated from Russian), {\it  Russian Acad.\
Sci.\ Sb.\ Math.} 82(1995), 1--20.\

\bibitem[BP4]{BP4} Bedford, E.\ and Pinchuk, S.,
  Domains in ${\CC}^2$ with non-compact automorphism group,
preprint.

\bibitem[Ber1]{Ber1}  Berteloot, F., Sur certains domaines faiblement
pseudoconvexes dont le groupe
d'automorphismes analytiques est non compact,
{\it  Bull.\ Sci.\ Math.} (2) 114(1990), 411--420.

\bibitem[Ber2]{Ber2} Berteloot, F., Un principe de localisation
pour les domaines faiblement
pseudoconvexes de ${\CC}^2$ dont le groupe d'automorphismes
holomorphes est non compact, {\it Colloque d'Analyse Complexe
et G\'eom\'etrie} (Marseille, 1992), {\it Asterisque} 217(1993),
13--27.

\bibitem[Ber3]{Ber3} Berteloot, F.,
Characterization of models in ${\CC}^2$
by their automorphism groups, {\it Internat.\ J.\ Math.} 5(1994),
619--634.\

\bibitem[BeCo]{BeCo} Berteloot, F. and C{\oe}ur\'e G., Domaines de ${\CC}^2$,
pseudoconvex et de type fini ayant un groupe non compact
d'automorphismes, {\it Ann.\ Inst.\ Fourier Grenoble} 41(1991), 77--88.

\bibitem[BM]{BoMa} Bochner, S. and Martin, W., Several complex variables,
{\it Princeton Univ.\ Press}, 1948.

\bibitem[FIK1]{FIK1} Fu, S., Isaev, A.\ V. and Krantz, S.\ G.,
  Examples of
domains with non-compact automorphism groups, {\it Math.\ Res.\ Letters}
3(1996), 609--617.\

\bibitem[FIK2]{FIK2}  Fu, S., Isaev, A.\ V. and
Krantz, S.\ G., Reinhardt domains
with non-compact automorphism groups,
{\it Math.\ Res.\ Letters} 3(1996), 109--122.

\bibitem[Fr]{Fr}  Frankel, S.,
Complex geometry of convex domains that cover varieties,
{\it Acta Math.} 163(1989), 109--149.\

\bibitem[GK1]{GK1}  R.\ E.\ Greene and S.\ G.\ Krantz, Characterizations of
certain weakly pseudo-convex domains with non-compact automorphism
groups, in {\it Complex Analysis Seminar}, Springer Lecture Notes
Vol.\ 1268 (1987), 121-157.

\bibitem[IK1]{IKR}  A.\ Isaev and S.\ G.\ Krantz, Domains with non-compact
automorphism group:  A survey, {\it Advances in Math.} 146(1999), 1--38.

\bibitem[IK2]{IK1}    Isaev, A.\ V.\ and Krantz,
S.\ G., On the boundary orbit accumulation set for
a domain with non-compact automorphism group,
 {\it Michigan Math.\ J.} 43(1996), 611--617.

\bibitem[IK3]{IK2}    Isaev, A.\ V.\ and Krantz,
S.\ G., Finitely smooth Reinhardt domains with non-compact
automorphism group, {\it Illinois.\ J.\ Math.} 41(1997), 412--420.

\bibitem[IK4]{IK3}  Isaev, A.\ V.\ and Krantz,
S.\ G., Hyperbolic Reinhardt domains with non-compact
automorphism group, {\it Pacific J.\ Math.} 184(1998), 149--160.


\bibitem[Ki1]{Ki1}   Kim, K.-T.,
  Domains with non-compact automorphism groups, {\it
Recent Developments in Geometry} (Los Angeles, CA, 1987), 249--262,
{\it Contemp.\ Math.} 101, Amer.\ Math.\ Soc.,1989.

\bibitem[Ki2]{Ki2}  Kim, K.-T.,
Complete localization of domains with
non-compact automorphism groups, {\it
Trans.\ Amer.\ Math.\ Soc.} 319(1990), 139--153.\

\bibitem[Ki3]{Ki3}   Kim, K.-T.,
Domains in ${\CC}^n$ with a piecewise Levi
flat boundary which possess a non-compact automorphism group,
{\it Math.\ Ann.} 292(1992), 575--586.\

\bibitem[Ki4]{Ki4}  Kim, K.-T.,
Geometry of bounded domains and the scaling techniques in
several complex variables, {\it Lecture Notes Series} 13, Seoul
National University, Research Institute of Mathematics, Global
Analysis Research Center, Seoul, 1993.

\bibitem[Ki5]{Ki5}   Kim, K.-T.,
  On a boundary point repelling
automorphism orbits, {\it J.\ Math.\ Anal.\ Appl.} 179(1993),
463--482.

\bibitem[Ki6]{Ki6}  Kim, K.-T.,
Two examples for scaling methods in several complex variables, {\it
RIM-GARC Preprint Series, Seoul National University} 95-53(1995).

\bibitem[KIKR]{KIKR}  Kim, K.-T. and Krantz, S.\ G.,
A crash course in the function
theory of several complex variables, {\it Complex geometric analysis in
Pohang}(1997), 3--37, {\it Contemp.\ Math.} 222, American Math.\ Society,
Providence, RI, 1999.

\bibitem[KRA1]{KRA1}  Krantz, S.\ G., {\it Function Theory of Several Complex
Variables}, $2^{\rm nd}$ Ed., Wadsworth, Belmont, 1992.

\bibitem[KRA2]{KRA2}  Krantz, S.\ G., Convexity in complex analysis,
{\it Proc.\ Symp.\ Pure Math.}, vol. 52 (E.\ Bedford, J.\ D'Angelo,
R.\ Greene, and S.\ Krantz eds.), American Mathematical Society,
Providence, R.I., 1991.

\bibitem[LEM1]{LEM1}  Lempert, L., La metrique de Kobayashi et la
representation des domains sur
la boule, {\it Bull.\ Soc.\ Math.\ France} 109(1981), 427--474.

\bibitem[LEM2]{LEM2}  Lempert, L., The Dolbeault complex in infinite
dimensions, {\it J.\ AMS} 11(1998), 485--520.

\bibitem[Mu]{Mu} Mujica, J., Complex analysis in Banach spaces,
North-Holland, 1986.

\bibitem[NAR]{NAR}  Narasimhan, R., {\it Several Complex Variables},
University of Chicago Press, Chicago, 1971.

\bibitem[ROS]{ROS}  Rosay, J.-P.,
Sur une characterization de la boule parmi les domains de
$\CC^n$ par son groupe d'automorphismes, {\it Ann.\ Inst.\ Four.\
Grenoble} XXIX(1979), 91--97.

\bibitem[WON]{WON}  Wong, B., Characterization of the ball in $\CC^n$ by
its automorphism group, {\it Invent.\ Math.} 41(1977), 253--257.
\end{thebibliography}
\end{document}